\documentclass[12pt, reqno]{amsart}

\usepackage{amsfonts}
\usepackage{amssymb}
\usepackage{latexsym}
\usepackage{epsf,graphicx}
\usepackage{amsmath}

\newtheorem{remar}{Remark}
\newtheorem{examp}{Example}

\newcommand{\rema}{\begin{remar}\rm}
\newcommand{\erema}{$\blacktriangleright$\end{remar}\vspace{2mm}}

\newcommand{\exa}{\begin{examp}\rm}
\newcommand{\eexa}{$\blacktriangleright$\end{examp}\vspace{2mm}}

\newcommand {\Hol}{\mathop{\rm Hol}\nolimits}

\renewcommand {\Im}{\mathop{\rm Im}\nolimits}
\renewcommand {\Re}{\mathop{\rm Re}\nolimits}
\newcommand{\pl}{\partial}
\newcommand{\pr}{\noindent{\bf Proof.}\quad }
\newcommand{\epr}{\ $\blacksquare$ \vspace{2mm} }

\newcommand{\be} {\begin{eqnarray}}
\newcommand{\ee} {\end{eqnarray}}
\newcommand{\bep} {\begin{eqnarray*}}
\newcommand{\eep} {\end{eqnarray*}}

\newfont{\bbb}{msbm10 at 12pt}
\def\Bbb#1{\hbox{\bbb #1}}

\newcommand{\C}{{\Bbb C}}

\newtheorem{defin}{Definition}%[section]
\newtheorem{theorem}{Theorem}%[section]
\newtheorem{corol}{Corollary}%[section]

\begin{document}

\title{Covering results and perturbed Roper--Suffridge operators}

\author[M. Elin]{Mark Elin}
\address{Department of Mathematics,
Ort Braude College, Karmiel 21982, Israel}
\email{mark-elin@braude.ac.il}

\author[M. Levenshtein]{Marina Levenshtein}
\address{Department of Mathematics,
Ort Braude College, Karmiel 21982, Israel}
\email{marlev@braude.ac.il}


\begin{abstract}
This work is devoted to the advanced study of Roper--Suffridge
type extension operators. For a given non-normalized spirallike
function (with respect to an interior or boundary point) on the
open unit disk of the complex plane, we construct perturbed
extension operators in a certain class of Banach spaces and prove
that these operators preserve the spirallikeness property. In
addition, we present an extension operator for semigroup
generators.

We use a new geometric approach based on the connection between
spirallike mappings and one-parameter continuous semigroups. It
turns out that the new one-dimensional covering results
established below are crucial for our investigation.
\end{abstract}

\maketitle

\section{Introduction}

In the last fifteen years, many authors have devoted their works
to the study of extension operators in the following sense. Given
a family of biholomorphic mappings on the unit ball of a proper
subspace of a Banach space $X$, to extend it to a family of
biholomorphic mappings of the unit ball of the whole space $X$
with preserving some required geometric properties such as
convexity, starlikeness, spirallikeness and others. These
investigations were initiated in 1995 with the paper by Roper and
Suffridge \cite{Ro-Su}, where they introduced an extension
operator which preserves a variety of required geometric
characteristics. Namely, given a (locally) univalent function $h$
on the open unit disk $\Delta :=\{z:|z|<1\}$ of the complex plane
$\mathbb{C}$ normalized by $h(0)=h'(0)-1=0$, they considered the
mapping $\Phi[h]:\mathbb{B}^n\mapsto\mathbb{C}^n$ on the open unit
ball $\mathbb{B}^{n}$ in $\mathbb{C}^{n}$ defined as follows:
\begin{equation}\label{RS}
h\mapsto \Phi[h](x,y)=\left(h(x),\sqrt{h'(x)}\,y \right),
\end{equation}
where $y\in\C^{n-1},\ |x|^2+\|y\|^2<1$. This operator preserves
convexity \cite{Ro-Su}, starlikeness \cite{G-K-2000},
$\mu$-spirallikeness \cite{G-K-K-2000, L-L-2005}, as well as other
important properties.

Later, this idea was developed in different directions (see, for
example, \cite{G-K-book, F-L, E-2011} and references therein),
diverse generalizations and modifications of this operator where
considered in various normed spaces.

Recently, the first author \cite{E-2011} introduced a general
construction which involves most extension operators considered
earlier, namely, the operator
\begin{equation}\label{E}
h\mapsto\Phi _{\Gamma} [h](x,y):=(h(x),\Gamma(h,x)y)
\end{equation}
defined on the unit ball of the direct product $X\times Y$ of two
Banach spaces, where $\Gamma (h,x)$ is an operator-valued mapping
satisfying some natural conditions. It was proved in \cite{E-2011}
that these conditions provide that the extension operator $\Phi
_{\Gamma}[h]$ preserves starlikeness and spirallikeness.

In \cite{M1, M2}, Muir considered an extension operator which is
not included in construction (\ref{E}). More precisely, given a
univalent function $h$ on $\Delta $ normalized by
$h(0)=h'(0)-1=0$, he introduced the operator
$\widehat{\Phi}[h]:\mathbb{B}^n\mapsto\mathbb{C}^n$ defined by
\begin{equation}\label{M}
h\mapsto
\widehat{\Phi}[h](x,y)=\left(h(x)+h'(x)Q(y),\sqrt{h'(x)}\,y
\right),
\end{equation}
where $Q$ is a homogeneous quadratic polynomial and proved, in
particular, that $\widehat{\Phi}$ preserves starlikeness whenever
$\displaystyle\sup_{\|y\|=1}|Q(y)|\leq \frac{1}{4}\,$. The proof
of the last result is carried out  in \cite{M1} by a technical
verification of Suffridge's starlikeness criterion (see
\cite{STJ-77}).

In Section \ref{geom}, we give a geometric explanation of this
fact. By the way, this answers question (c) in \cite{E-2011}. It
turns out that the mentioned above Muir's theorem follows
immediately from some new one-dimensional covering theorems
established in Section \ref{cover}. Section \ref{Roper} is devoted
to an advanced study of the Roper--Suffridge type extension
operators. Moreover, the semigroup approach enables us to study
special perturbed extension operators in a certain class of Banach
spaces in Section \ref{muir}. These operators preserve
spirallikeness (with respect to an arbitrary interior or boundary
point).

\section{Preliminary notions}

In this section we present some notions of nonlinear analysis and
geometric function theory which will subsequently be useful.

Let $X$ be a complex Banach space with the norm $\|\cdot\|$.
Denote by $\Hol(D,E)$ the set of all holomorphic mappings on a
domain $D\subseteq X$ which map $D$ into a set $E\subseteq X$, and
set $\Hol(D):=\Hol(D,D)$.

We start with the notion of a one-parameter continuous semigroup.

\begin{defin}\label{def-sg}
Let $D$ be a domain in a complex Banach space $X$. A family
$S=\{F_t\}_{t\ge0}\subset\Hol(D)$ of holomorphic self-mappings of
$D$ is said to be a one-parameter continuous semigroup (in short,
semigroup) on $D$ if \be\label{semi_prop1} F_{t+s}=F_t\circ
F_s,\quad t,s\geq 0, \ee and for all $x\in D$,
\begin{eqnarray}\label{semi_prop2}        %equa1.3
\lim_{t\rightarrow 0^+}F_t(x) =x.
\end{eqnarray}
\end{defin}

\begin{defin} \label{def-gen}
A semigroup $S=\{F_t\}_{t\ge0}$ on $D$ is said to be generated if
for each $x\in D$, there exists the strong limit
\begin{eqnarray}
f(x):=\lim_{t\rightarrow 0^{+}}\,\frac 1t\biggl(x-F_t(x)\biggr).
\end{eqnarray}

In this case the mapping $f:\ D\mapsto X$ is called the
(infinitesimal) generator of the semigroup $S$.
\end{defin}

In fact, semigroups are one of the key tools in the geometric
function theory. For instance, a mapping $h$ is spirallike if
there exists a linear semigroup on its image. More precisely:

\begin{defin}%[see \cite{E-R-S-2004, RS-SD-book}, cf., \cite{STJ-77, GS-99, G-K-book}]\label{def-spiral}
Let $h$ be a biholomorphic mapping defined on a domain $D$ of a
Banach space $X$. The mapping $h$ is said to be spirallike if
there is a bounded linear operator $A$ such that the function
$\Re\lambda$ is bounded away from zero on the spectrum of $A$ and
such that for each point $w\in h(D)$ and each $t\ge0$, the point
$e^{-tA}w$ also belongs to $h(D)$. In this case $h$ is called
$A$-spirallike. If $A$ can be chosen to be the identity mapping,
that is, $e^{-t}w\in h(D)$ for all $w\in h(D)$ and all $t\ge0$,
then $h$ is called starlike.
\end{defin}

It follows by this definition that the origin belongs to the
closure $\overline{ h(D)}$. In the case where the image $h(D)$
contains the origin, we say that $h$ is spirallike with respect to
an interior point. Otherwise, if the origin belongs to the
boundary of $h(D)$, $h$ is said to be spirallike with respect to a
boundary point.

In the one-dimensional case, a univalent function
$h:\Delta\mapsto\mathbb{C}$ is $\mu$-spirallike ($\Re \mu >0$) if
for all $t\geq 0$, $$e^{-\mu t}h(\Delta)\subseteq h(\Delta).$$ If
this inclusion holds for $\mu =1$, the function $h$ is starlike.

\section{Geometric interpretation} \label{geom}

Let $h$ be a normalized starlike function on the open unit disk of
the complex plane. By Theorem 4.1 in \cite{M1}, the mapping
$\widehat{\Phi}[h]$ defined by~(\ref{M}) is the starlike mapping
of the open unit ball in $\C^n$ whenever
$\displaystyle\sup_{\|y\|=1}|Q(y)|\leq \frac{1}{4}$\,. This means
that its image is invariant under the action of the semigroup
$\{G_t\}_{t\ge0}$, where
\[
G_t(z,w)=\left(e^{-t}z,e^{-t}w  \right),\quad z\in\C,\
w\in\C^{n-1}.
\]

In addition, it is easy to see that
$\widehat{\Phi}[h]=\Psi\circ\Phi[h]$, where $\Phi[h]$ is the
Roper--Suffridge extension operator defined by~(\ref{RS}) and
$\Psi(z,w)=\left(z+Q(w),w\right)$ is the automorphism of the space
$\mathbb{C}\times\mathbb{C}^{n-1}$. Therefore, the family
$\{F_t\}_{t\ge0}$ defined by $F_{t}:=\Psi^{-1}\circ G_t\circ\Psi$
forms a semigroup on the image of the mapping $\Phi[h]$. A direct
calculation shows that
\begin{equation}\label{sg}
F_t(z,w)=\left(e^{-t}z + \left(e^{-t}-e^{-2t}\right)Q(w),e^{-t}w
\right).
\end{equation}

On the other hand, it is known  that the mapping $\Phi[h]$ is
starlike \cite{G-K-2000}, i.e., its image is invariant under the
action of the semigroup
\begin{equation}\label{star}
(z,w)\mapsto\left(e^{-t}z,e^{-t}w \right).
\end{equation}

Comparing semigroup elements (\ref{sg}) and (\ref{star}), we see
that both of them have the same second coordinate.

Fix a point $(z_0,w_0)$ in the image of $\Phi [h]$ and $t>0$.
Consider the preimage of the intersection of
$\Phi[h](\mathbb{B}^n)$ with $w=e^{-t}w_0$. This preimage is a
submanifold of the ball $\mathbb{B}^n$. Denote by $\Omega$ the
orthogonal projection of this submanifold on the first coordinate
plane. Taking into account that the first coordinate of the
Roper--Suffridge extension operator is $h(x)$, we conclude that
the mentioned Muir's result is equivalent to the fact that
$h(\Omega)$ covers some disk centered at~$e^{-t}z_0$.

According to this explanation,  we start our consideration from
two auxiliary one-dimensional covering results which are of
independent interest.

\section{Covering results}\label{cover}

In this section, we consider non-normalized univalent functions
$h\in\Hol(\Delta,\mathbb{C})$ on the open unit disk $\Delta$. For
some special subsets $\Omega _{\alpha}$ of $\Delta$, we determine
disks covered by $h(\Omega _{\alpha})$.

\begin{theorem}\label{t1}
Let $h\in \mathrm{Univ} (\Delta , \mathbb{C})$. For $\alpha\in
(0,1)$, and $x_{0}\in \Delta$ define the set
\[
\Omega _{\alpha} :=\left\{ x\in\Delta :\alpha
|h'(x_{0})|(1-|x_{0}|^{2})<|h'(x)|(1-|x|^{2}) \right\}.
\]
Then the image $h(\Omega _{\alpha})$ covers the open disk of
radius $\displaystyle\frac{1-\alpha
}{4}|h'(x_{0})|(1-|x_{0}|^{2})$ centered at $h(x_{0})$.
\end{theorem}

\pr Consider the function $g:\Delta\mapsto\mathbb{C}$ defined by
\[
g(z):=\frac{h\!\!\left(\displaystyle\frac{x_{0}-z}{1-\overline{x_{0}}z}\right)-h(x_{0})}{h'(x_{0})(|x_{0}|^{2}-1)},
\quad z\in\Delta .
\]
Since $g$ is a normalized univalent function, by the Koebe
One-Quarter Theorem, $g(\Delta)$ covers the disk $|z|<\frac{1}{4}$
and, consequently, $h(\Delta)$ covers the open disk of radius
$\frac{1}{4}|h'(x_{0})|(1-|x_{0}|^{2})$ centered at $h(x_{0})$.

Standard calculations show that
$x:=\displaystyle\frac{x_{0}-z}{1-\overline{x_{0}}z}\in\Omega_\alpha$
if and only if
\begin{equation} \label{a1}
z\in\widetilde{\Omega}_\alpha:=\left\{ z\in\Delta :
|g'(z)|>\frac{\alpha }{1-|z|^{2}}\right\}.
\end{equation}

On the other hand, by the Koebe Distortion Theorem, for each
$z\in\Delta$,
\begin{equation} \label{a2}
|g'(z)|\geq \frac{1-|z|}{(1+|z|)^{3}}
\end{equation}
and
\begin{equation} \label{a3}
|g(z)|\geq \frac{|z|}{(1+|z|)^{2}}.
\end{equation}
It follows from (\ref{a1}) and (\ref{a2}), that for each
$z\in\Delta \setminus \widetilde{\Omega}_\alpha$,
\[
\frac{1-|z|}{(1+|z|)^{3}}\leq |g'(z)|\leq \frac{\alpha
}{1-|z|^{2}},
\]
hence,
\[
(1-\alpha)|z|^{2}-2(1+\alpha)|z|+(1-\alpha)\leq 0 ,
\]
which implies
$|z|\ge\displaystyle\frac{1-\sqrt{\alpha}}{1+\sqrt{\alpha}}$.

Therefore, for each point $z\in\Delta\setminus
\widetilde{\Omega}_\alpha$, by inequality (\ref{a3}),
\[
|g(z)|\geq \frac{|z|}{(1+|z|)^{2}}\geq
\frac{\frac{1-\sqrt{\alpha}}{1+\sqrt{\alpha}}}{\left(1+\frac{1-\sqrt{\alpha}}{1+\sqrt{\alpha}}\right)^{2}}=\frac{1-\alpha
}{4}.
\]

Let now $x^{*}\in\Delta\setminus\Omega _{\alpha}$. Then
$z^{*}:=\displaystyle\frac{x_{0}-x^{*}}{1-\overline{x_{0}}x^{*}}\in\Delta\setminus
\widetilde{\Omega}_\alpha$ and hence
$$|g(z^{*})|=\frac{|h(x^{*})-h(x_{0})|}{|h'(x_{0})|(1-|x_{0}|^{2})}\geq
\frac{1-\alpha }{4} , $$ or
$$|h(x^{*})-h(x_{0})|\geq\frac{1-\alpha
}{4}|h'(x_{0})|(1-|x_{0}|^{2}),$$ and the assertion follows. \epr

\begin{theorem}\label{t2}
Let $\alpha\in\mathbb{R}$ and $\beta\in\mathbb{C}$ be such that
$0<\alpha<|\beta |<1$, and suppose that $h\in \mathrm{Univ}
(\Delta , \mathbb{C})$ satisfies $\beta h(\Delta)\subseteq
h(\Delta)$. Given a point $x_{0}\in \Delta$, define the point
$x_1=h^{-1}\left(\beta h(x_0) \right)\in\Delta$ and the set
\[
\Omega_\alpha := \left\{ x \in \Delta : \alpha
|h'(x_{0})|(1-|x_{0}|^{2}) < |h'(x)|(1-|x|^{2})\right\} .
\]
Then the image $h(\Omega_{\alpha})$ covers the disk centered at
$h(x_1)=\beta h(x_{0})$ of radius  %\linebreak
$\displaystyle\frac{|\beta|-\alpha}{4|\beta|}|h'(x_{1})|(1-|x_{1}|^{2})$
which is not less than
$R=\displaystyle\frac{|\beta|-\alpha}{4}|h'(x_0)|(1-|x_0|^2)$.
\end{theorem}

\pr By our assumptions, the function $F:=h^{-1}(\beta h)$ is a
holomorphic self-mapping of $\Delta$ and, consequently, satisfies
the Schwarz--Pick inequality:
\[
|F'(x_{0}| \leq \frac{1-|F(x_{0})|^{2}}{1-|x_{0}|^{2}}=
\frac{1-|x_1|^{2}}{1-|x_{0}|^{2}}  \,.
\]
Set $g:=h \circ F =\beta h$. Then
\begin{equation}\label{gtag}
|g'(x_{0})|(1-|x_{0}|^{2}) = |h'(x_{1})|\cdot
|F'(x_{0})|(1-|x_{0}|^{2}) \leq |h'(x_{1})|(1-|x_{1}|^{2}) .
\end{equation}
On the other hand, $|g'(x_{0})|=|\beta|\cdot |h'(x_{0})|$ and so
\[
|\beta|\cdot |h'(x_{0})|(1-|x_{0}|^{2})\leq
|h'(x_{1})|(1-|x_{1}|^{2}) .
\]
Hence, the set
\[
\Omega := \left\{ x\in\Delta: \frac{\alpha}{|\beta|}\cdot
|h'(x_{1})| (1-|x_{1}|^{2}) < |h'(x)| (1-|x|^{2}) \right \}
\]
is a subset of $\Omega _{\alpha}$. By Theorem \ref{t1}, the disk
centered at $h(x_{1})$ of radius $\displaystyle
\frac{1-\frac{\alpha}{|\beta|}}{4}
\left|h'(x_{1})\right|(1-|x_{1}|^{2})$ lies in $h\left(\Omega
\right)$, and the assumption follows. \epr

\rema Note that the condition $\beta h(\Delta)\subseteq h(\Delta)$
with $|\beta |<1$ is necessary and sufficient for a holomorphic
mapping $h:\Delta \mapsto \mathbb{C}$ to be a solution to the
so-called Schr\"oder functional equation
\begin{equation}\label{sch}
h\circ F=\beta h,
\end{equation}
where $F$ is a holomorphic self-mapping of the open unit disk
$\Delta$ of dilation or hyperbolic type. In the dilation case, the
Denjoy--Wolff point $\tau$ of $F$ is an interior point of
$\Delta$. It is clear that if (\ref{sch}) has a solution, then
$\beta=F'(\tau)$. It was proved by K\oe nigs that the last
equality is also sufficient for the solvability of the Schr\"oder
equation~(\ref{sch}) whenever $|\beta|\not=0,1$. Otherwise, when
$F$ is of hyperbolic type, the Denjoy--Wolff point $\tau$ belongs
to the boundary $\pl{\Delta}$ and
$F'(\tau):=\displaystyle\lim_{r\to1^-}
\frac{F(r\tau)-\tau}{\tau(r-1)}$ exists and is a positive number
less than $1$. For the case $\beta=F'(\tau)$ the existence of a
solution to equation (\ref{sch}) was proved by Valiron. The survey
of the related topics can be found in \cite{E-S-book}.

On the other hand, if the function $h$ is $\mu$-spirallike then
the inclusion $\beta h(\Delta)\subseteq h(\Delta )$ holds for each
$\beta$ of the form $\beta=e^{-t\mu},\ t>0$. \erema

To finish this section, we remark that according to the geometric
explanation given in Section~\ref{geom}, Muir's theorem on
starlikeness preserving follows from Theorem~\ref{t2} with
$\alpha=e^{-2t}$ and $\beta=e^{-t}$. In the next section, we
consider a more general situation.

\section{Roper--Suffridge type extensions}\label{Roper}

Let $(Y,\|\cdot\|_{Y})$ be a complex Banach space and let $r\ge1$.
Consider the set
\[
\mathbb{B}:=\left\{ (x,y)\in\mathbb{C}\times Y:
|x|^{2}+\|y\|_{Y}^{r}<1\right\}.
\]
This set is the open unit ball of $\mathbb{C}\times Y$ with
respect to a norm $\|\cdot\|$. Actually, this norm is the
Minkowski functional of the set $\mathbb{B}$. Obviously,
$\mathbb{C}\times Y$ equipped with this norm $\|\cdot\|$ is a
complex Banach space.

It follows from Theorem 5.1 in \cite{E-2011} that the extension
$H(x,y)=(h(x),h'(x)^{\frac{1}{r}}y)$ of a $\mu$-spirallike
function $h$ on $\Delta $ is $\left (\begin{smallmatrix} \mu & 0
\\ 0 & B+\frac{\mu}{r}\mathrm{id}_{Y}
\end{smallmatrix}\right )$-spirallike for each linear accretive operator $B$ on $Y$
(note in passing that the branch of the power
$h'(x)^{\frac{1}{r}}$ can be chosen appropriately). We apply
covering results from the previous section to the advanced study
of this operator. Our main result (Theorem~\ref{main}) provides a
new tool for the study of various extension operators. In
particular, in the next section we apply it to prove that Muir's
operator preserves spirallikeness in a certain class of Banach
spaces.
%Extension operators in Corollaries~\ref{cor1}
%and~\ref{cor2} can be understood as perturbations of the
%Roper--Suffridge operator preserving spirallikeness.

\begin{theorem}\label{main}
Let $h:\Delta\mapsto\mathbb{C}$ be a $\mu$-spirallike mapping and
$H:\mathbb{B}\mapsto\mathbb{C}\times Y$ be defined by
\begin{equation}\label{r}
H(x,y):=\left(h(x), h'(x)^{\frac{1}{r}}y \right) .
\end{equation}
Suppose that $B\in L(Y)$ generates a semigroup of strict
contractions. Then for each point $(z_{0},w_{0})\in H(\mathbb{B})$
and for all $t>0$,
\[
\left( e^{-\mu t}z_{0}+\gamma _{t}, e^{-\frac{\mu}{r}t}e^{-Bt}
w_{0}\right) \in H(\mathbb{B})
\]
whenever $|\gamma _{t}
|<R_{t}=\frac{1-\|e^{-Bt}\|^r}{4}|h'(x_{1})|(1-|x_{1}|^{2})$ with
$x_{1}=h^{-1}(e^{-\mu t}z_{0})$.
\end{theorem}

Note that if the operator $B$ is scalar, i.e., $B=\lambda\,{\rm
id}_Y$ ($\lambda\in\mathbb{C}$), it generates a semigroup of
strict contractions if and only if $\Re\lambda>0$.

\pr Since $H$ is a $\left (\begin{smallmatrix} \mu & 0
\\ 0 & B+\frac{\mu}{r}\mathrm{id}_{Y}
\end{smallmatrix}\right )$-spirallike mapping (see \cite{E-2011}),
its image $H(\mathbb{B})$ is invariant under the action of the
semigroup $\{F_{t}\}_{t \geq 0}$ defined by
\begin{equation}\label{ftequation}
F_{t} (z,w) = \left(e^{-\mu t} z, e^{-\frac{\mu}{r}t}e^{-Bt}
w\right) , \quad (z,w) \in H(\mathbb{B}), \quad t \geq 0.
\end{equation}
Fix $t>0$ and $(z_{0}, w_{0}) \in H(\mathbb{B})$. Then $z_{0} =
h(x_{0})$ and $w_{0} = h'(x_{0})^{\frac{1}{r}} y_{0}$ for some
$(x_{0}, y_{0}) \in \mathbb{B}$.

Denote $(z_{1},w_{1}):=\left(e^{-\mu
t}z_{0},e^{-\frac{\mu}{r}t}e^{-Bt}w_{0}\right) \in H(\mathbb{B})$.
Then there is a point $(x_{1},y_{1}) \in \mathbb{B}$ such that
$z_{1}=h(x_{1})$ and $w_{1} = h'(x_{1})^{\frac{1}{r}} y_{1}$. Let
$U$ be a region in the complex plane $\mathbb{C}$ containing
$z_{1}$ and such that for all $z\in U$, the point $(z,w_{1})$
belongs to the image $H(\mathbb{B})$. We have to estimate disks
contained in~$U$.

Consider the set of all $(x,y) \in \mathbb{B}$ such that
$h'(x)^{\frac{1}{r}} y=w_1$ and, consequently, $h(x)\in U$. Since
$|x|^{2}+\|y\|_{Y}^{r} < 1$, we have $\|w_{1}\| ^{r} <
|h'(x)|(1-|x|^{2}).$

On the other hand,  $w_{1} = e^{-\frac{\mu}{r}t}e^{-Bt} w_{0} =
e^{-\frac{\mu}{r}t}e^{-Bt} h'(x_{0})^{\frac{1}{r}} y_{0}$.
Therefore,
\begin{equation}\label{inequality}
|e^{-\mu t}|\cdot\|e^{-Bt}\|^{r}\cdot |h'(x_{0})| \cdot \| y_{0}
\|_{Y} ^{r} < |h'(x)|(1-|x|^{2}).
\end{equation}
Define the set
\[%begin{equation}\label{omegaequation}
\Omega := \left \{ x \in \Delta : |e^{- \mu
t}|\cdot\|e^{-Bt}\|^{r}\cdot |h'(x_{0})| \cdot (1-|x_{0}|^{2}) <
|h'(x)|(1-|x|^{2}) \right \}.
\]%end{equation}

It is clear that for each $x \in \Omega$ and
$y=\displaystyle\frac{1}{h'(x)^{\frac{1}{r}}}\,w_{1}$, we have
$(x,y)\in \mathbb{B}$ and $h'(x)^{\frac{1}{r}} y = w_{1}$. Thus,
$h(\Omega)\subseteq U$.

By Theorem \ref{t2} with $\beta =e^{-\mu t}$ and $\alpha =|e^{-\mu
t}|\cdot\|e^{-Bt}\|^{r}$, the image $h(\Omega)$ covers the open
disk centered at $z_{1} = e^{-\mu t} z_{0} $ of radius
$\displaystyle \frac{1- \|e^{-Bt}\|^{r}}{4} |h'(x_{1})| \cdot
(1-|x_{1}| ^{2})$. The proof is complete.   \epr

\begin{corol}\label{cor1}
Let $h:\Delta\mapsto\mathbb{C}$ be $\mu$-spirallike and $A$ be a
linear operator in $\mathbb{C}\times Y$ with $\Re
\sigma>\varepsilon >0$ for each $\sigma$ from the spectrum of $A$.
Suppose that  there exists an automorphism $\Phi$ of the space
$\mathbb{C}\times Y$ such that for each point $(z_{0},w_{0})\in
H(\mathbb{B})$ and for all $t\geq 0$,
\begin{equation}\label{c1}
\Phi^{-1}\left(e^{-tA}\Phi
(z_{0},w_{0})\right)=\left(z_{t},e^{-\frac{\mu}{r}t}e^{-Bt}w_{0}\right),
\end{equation}
for some $B\in L(Y)$ which generates a semigroup of strict
contractions. If
\begin{equation}\label{c2} \left|z_{t}-e^{-\mu
t}z_{0}\right|\leq
R_{t}=\frac{1-\|e^{-Bt}\|^{r}}{4}|h'(x_{1})|(1-|x_{1}|^{2}),
\end{equation}
where $x_{1}=h^{-1}(e^{-\mu t}z_{0})$, then the mapping $\Phi\circ
H: \mathbb{B}\mapsto \mathbb{C}\times Y$ is $A$-spirallike.
\end{corol}

\pr Denote $\gamma _{t}:=z_{t}-e^{-\mu t}z_0$ and rewrite equality
(\ref{c1}) in the form
\[
\Phi^{-1}\left(e^{-tA}\Phi (z_{0},w_{0})\right)=\left(e^{-\mu t}
z_0 +\gamma_{t}, e^{-\frac{\mu}{r}t}e^{-Bt}w_{0}\right).
\]
Since $|\gamma _{t}|\leq R_{t}$ by (\ref{c2}), it follows from
Theorem \ref{main} that for all $t\geq 0$,
\[
\Phi^{-1}\left(e^{-tA}\Phi (z_{0},w_{0})\right)\subseteq
H(\mathbb{B})
\]
or, which is the same, $e^{-tA}\Phi (z_{0},w_{0}) \subseteq \Phi
(H(B)).$ \epr

\rema  It follows from the proof of Theorem~\ref{t2} that $R_t$
defined in Theorem~\ref{main} and Corollary~\ref{cor1} is not less
than $|e ^{-\mu t}|\frac{1- \|e^{-Bt}\|^r}{4} |h'(x_{0})| \cdot
(1-|x_{0}| ^{2})$. \erema

\section{Muir's type perturbations of extension
operators}\label{muir}

In this section, we present some generalizations of Muir's
perturbations of the Roper--Suffridge operators preserving
spirallikeness. It is worth mentioning that in saying
`spirallike', we mean spirallikeness with respect to an arbitrary
interior or boundary point and do not separate these two different
cases. Further, we construct extensions for holomorphic generators
on $\Delta$ of hyperbolic and dilation type which, in their turn,
generate semigroups on the unit balls of the corresponding Banach
spaces.

As above, $(Y, \|\cdot\|_{Y})$ is a complex Banach space, and
$\mathbb{C}\times Y$ is the Banach space with the open unit ball
$\mathbb{B}=\{(x,y)\in\mathbb{C}\times Y: |x|^{2}+\|y\|^{r}<1 \}$
for some $r\geq 1$. For simplicity, we will use the results of the
previous section in the particular case when $B$ is a scalar
operator $B=\lambda\,{\rm id}_Y$.

\begin{theorem}\label{cor2}
Let $h:\Delta\mapsto\mathbb{C}$ be $\mu$-spirallike and $Q: Y
\mapsto\mathbb{C}$ be a homogeneous polynomial of degree
$r\in\mathbb{N}$. Then the mapping
$\widehat{H}:\mathbb{B}\mapsto\mathbb{C}\times Y$ defined by
$$\widehat{H}(x,y):=\left(h(x)+h'(x)Q(y),h'(x)^{\frac{1}{r}}
y\right) $$ is $\displaystyle\left (\begin{smallmatrix} \mu & 0
\\ 0 & \left(\lambda +\frac{\mu }{r}\right){\rm id}_{Y}
\end{smallmatrix}\right )$-spirallike for each
$\lambda\in\mathbb{C},\ \Re \lambda
>0$, whenever
\[
\displaystyle\sup_{\|y\|=1}|Q(y)|\leq \frac{1}{4}\cdot\frac{\Re
\lambda }{|\lambda |}\,.
\]
Moreover, this bound is sharp.
\end{theorem}

\pr The mapping $\widehat{H}$ can be represented in the form
\[
\widehat{H}=\Phi \circ H ,
\]
where $\Phi (z,w):=(z+Q(w),w)$ is the automorphism of the space
$\mathbb{C}\times Y$ and $H$ is defined by (\ref{r}).

We have to show that for each element $F_{t}$ of the semigroup
$\{F_{t}\}_{t\geq 0}$ defined by (\ref{ftequation}), the inclusion
\[
F_{t}\circ\Phi (H(\mathbb{B}))\subseteq \Phi (H(\mathbb{B}))
\]
holds or, which is the same, for each point $(z_{0},w_{0})\in
H(\mathbb{B})$ and for all $t>0$,
\begin{eqnarray*}
&&\Phi ^{-1}\left(F_{t}\circ\Phi (z_{0},w_{0})\right)\\
&&=\left(e^{-\mu t}z_{0}+(e^{-\mu t}-e^{-\left(\lambda
+\frac{\mu}{r}\right) rt })Q(w_{0}),e^{-\left(\lambda
+\frac{\mu}{r}\right) t}w_{0}\right)\in H(\mathbb{B}).
\end{eqnarray*}

Fix $t>0$ and $(z_{0},w_{0})\in H(\mathbb{B})$, and denote $\gamma
_{t}:=(e^{-\mu t}-e^{-\left(\lambda +\frac{\mu}{r}\right) rt
})Q(w_{0})$. By our assumption, $$|Q(w_{0})|< \frac{1}{4}\frac{\Re
\lambda }{|\lambda |} |h'(x_{0})| \left(1-|x_{0}|^{2}\right),$$
where $x_{0}=h^{-1}(z_{0})$. Therefore, by Corollary \ref{cor1}
and Remark 2, it suffices to show that
\begin{equation}\label{cor}
\left|1-e^{-r\lambda t}\right|\frac{\Re \lambda}{|\lambda
|}<1-\left|e^{-r\lambda t}\right|.
\end{equation}

To this end, we denote
\[
f(t):=\left(\frac{1-|e^{-r\lambda t}|}{|1-e^{-r\lambda t
}|}\right)^{2}
\]
and prove that $\displaystyle\inf_{t>0}f(t)= \left(\frac{\Re
\lambda }{|\lambda|}\right)^{2}$.

Indeed, $\lim\limits _{t\rightarrow\infty}f(t)=1$ and $\lim\limits
_{t\rightarrow 0}f(t)=\left(\frac{\Re
\lambda}{|\lambda|}\right)^{2}$.

Moreover, if we denote $a:=\Re \lambda$ and $b:=\Im \lambda$, then
the derivative $f'(t)$ vanishes either when $\cos (brt)=1$ and
then $f(t)=1$, or when $\cos (brt)=\displaystyle
\frac{a^{2}(1+e^{-art})^{2}-b^{2}(1-e^{-art})^{2}}{a^{2}
(1+e^{-art})^{2}+b^{2}(1-e^{-art})^{2}}$ and then
\[
f(t)=\frac{a^{2}}{a^{2}+b^{2}}+\frac{b^{2}(1-e^{-art})^{2}}{a^{2}(1+e^{-art})
^{2}+b^{2}(1-e^{-art})^{2}}>\frac{a^{2}}{a^{2}+b^{2}}=\left(\frac{\Re
\lambda }{|\lambda|}\right)^{2}.
\]
The proof is complete. \epr

Now we consider extension operators for semigroup generators.

\begin{theorem} \label{gent}
Let $f: \Delta\mapsto\mathbb{C}$ be the generator of a semigroup
on $\Delta$ of dilation or hyperbolic type having the
Denjoy--Wolff point $\tau\in \overline\Delta$ and $Q\colon Y
\mapsto\mathbb{C}$ be a homogeneous polynomial of degree
$r\in\mathbb{N}$. Denote $\mu :=f'(\tau)$. Then for each
$\lambda\in\mathbb{C},\ \Re \lambda
>0$, the mapping $\widehat{f}:\mathbb{B}\mapsto\mathbb{C}\times
Y$ defined by
\begin{equation}\label{gen}
\widehat{f}(x,y):=\left(f(x)+Q(y),\frac{1}{r}\left(f'(x)+r\lambda
-\frac{\mu -f'(x)}{f(x)}Q(y) \right)y\right)
\end{equation}
is a semigroup generator on $\mathbb{B}$ whenever
$\sup\limits_{\|y\|=1}|Q(y)|\leq \displaystyle\frac{r\Re\lambda
}{4}$.
\end{theorem}

\pr Define the function $h\colon\Delta\mapsto\mathbb{C}$ as a
solution of the differential equation (see \cite{E-S-book})
\begin{equation} \label{h}
h'(z)f(z)=\mu h(z), \quad z\in\Delta .
\end{equation}
It is known that $h$ is a univalent function. Moreover, if $f$ is
the generator of a dilation type semigroup, the function $h$ is
$\mu$-spirallike with respect to an interior point. Otherwise, if
the generated semigroup is of hyperbolic type, $h$ is
$\mu$-spirallike with respect to a boundary point.

By Corollary \ref{cor2}, the mapping
$\widetilde{H}:\mathbb{B}\mapsto\mathbb{C}\times Y$ defined by
\[
\widetilde{H}(x,y):=\left(h(x)-\frac{h'(x)}{r\lambda }Q(y),
h'(x)^{\frac{1}{r}} y\right)
\]
is $\left (\begin{smallmatrix} \mu & 0
\\ 0 & \left(\lambda +\frac{\mu }{r}\right)\mathrm{id}_{Y}
\end{smallmatrix}\right )$-spirallike. Consequently, the family
$\{\widetilde{F}_{t}\}_{t\geq 0}$ of holomorphic mappings
$$\widetilde{F}_{t}(z,w):=\left(e^{-\mu t}z,e^{-t\left(\lambda
+\frac{\mu }{r}\right)id_{Y} }w\right),$$ is a semigroup on
$\widetilde{H}(\mathbb{B})$. Differentiating
$\widetilde{F}_{t}(z,w)$ at $t=0^{+}$, we find its generator
$$\widetilde{f}(z,w)=-\frac{d}{dt}\widetilde{F}_{t}(z,w)|_{t=0^{+}}=\left(\mu
z,\left(\lambda +\frac{\mu }{r}\right)w \right).$$

By Lemma 3.7.1 on p. 30 of \cite{E-R-S-2004}, the mapping
$\widehat{f}:\mathbb{B}\mapsto \mathbb{C}\times Y$ defined by
$$\widehat{f}(x,y):=[D\widetilde{H}(x,y)]^{-1}\widetilde{f}(\widetilde{H}(x,y)),$$
is a semigroup generator on $\mathbb{B}$.

By straightforward calculations we get
$$[D\widetilde{H}(x,y)]^{-1}=\left (\begin{smallmatrix}
\displaystyle\frac{1}{h'(x)} & \displaystyle\frac{1}{r\lambda
h'(x)^{\frac{1}{r}}}Q'(y)
\\ -\displaystyle\frac{h''(x)}{rh'(x)^{2}} & \;\displaystyle\frac{r\lambda h'(x)-h''(x)Q(y)}{r\lambda
h'(x)^{1+\frac{1}{r}}}id_{Y}
\end{smallmatrix}\right ),$$
and thus
\begin{equation}
\begin{split}
\widehat{f}(x,y)&=\left (\begin{smallmatrix}
\displaystyle\frac{1}{h'(x)} & \displaystyle\frac{1}{r\lambda
h'(x)^{\frac{1}{r}}}Q'(y)
\\ -\displaystyle\frac{h''(x)}{rh'(x)^{2}} & \;\displaystyle\frac{r\lambda h'(x)-h''(x)Q(y)}{r\lambda
h'(x)^{1+\frac{1}{r}}}id_{Y}
\end{smallmatrix}\right )\widetilde{f}(\widetilde{H}(x,y)) \\
&= \left(f(x)+Q(y),\frac{1}{r}(f'(x)+r\lambda -\frac{\mu
-f'(x)}{f(x)}Q(y) )y\right),
\end{split}
\end{equation}
which complete the proof. \epr

\rema Note that the condition $\mu =f'(\tau)$ was used in the proof of
Theorem~\ref{gent} %we only use
for the solvability of equation (\ref{h}) only. Since in the
hyperbolic case this equation is solvable for each $\mu $ such
that $|\mu -f'(\tau )|\leq f'(\tau ),\ \mu\not=0$ (see
\cite{E-S-book}), we conclude that in such a situation, our
assertion also holds. \erema

\rema In the last two theorems, we required for $r$ to be integer
by the only reason in order to define homogeneous polynomial $Q$
of degree $r$. At the same time, the polynomial $Q\equiv 0$
obviously satisfies all the conditions of the theorems. So, the
conclusions of Theorems~\ref{cor2} and~\ref{gent} hold with
$Q\equiv0$ and arbitrary $r\geq 1$. \erema

\end{document}